\documentclass[11pt]{article}
\usepackage{amsmath, amssymb, amscd, epsfig}

\newtheorem{thm}{Theorem}

\newtheorem{lm}[thm]{Lemma}

\newtheorem{rem}[thm]{Remark}

\newtheorem{ex}[thm]{Example}

\renewcommand{\phi}{\varphi}
\renewcommand{\epsilon}{\varepsilon}

\newcommand{\BB}{\mathbb}
\newcommand{\g}{\mathfrak}
\newcommand{\pf}{\noindent {\it Proof. }}
\newcommand{\qed}{\nopagebreak $\qquad$ $\square$ \vskip5pt}
\newcommand{\separate}{\vskip5pt}

\begin{document}

\title{\bf A New Proof of the Integral Localization Formula for
Equivariantly Closed Differential Forms}
\author{Matvei Libine}
\maketitle

\begin{abstract}
In this article we give a totally new proof of the
integral localization formula for equivariantly closed differential forms
(Theorem 7.11 in \cite{BGV}).
We restate it here as Theorem \ref{BGV}.

This localization formula is very well known, but the author hopes
to adapt this proof to obtain a more general result in the future.
\end{abstract}

\tableofcontents

\separate

\begin{section}
{Introduction}
\end{section}

Equivariant forms were introduced in 1950 by Henri Cartan.
There are many good texts on this subject including \cite{BGV}
and \cite{GS}.

Let $G$ be a compact Lie group acting on a compact oriented manifold $M$,
let $\g g$ be the Lie algebra of $G$,
and let $\alpha : \g g \to {\cal A}(M)$ be an equivariantly closed form on $M$.
For $X \in \g g$, we denote by $M_0(X)$ the set of zeroes of the vector
field on $M$ induced by the infinitesimal action of $X$.
We assume that $M_0(X)$ is discrete.
Then Theorem 7.11 in \cite{BGV} says
that the integral of $\alpha(X)$ can be expressed as a sum over the
set of zeroes $M_0(X)$ of certain {\em local} quantities of $M$ and $\alpha$:
$$
\int_M \alpha(X) = \sum_{p \in M_0(X)}
\text{local invariant of $M$ and $\alpha$ at $p$}.
$$
This is the essence of the integral localization formula for equivariantly
closed differential forms (Theorem \ref{BGV}).

\separate

This localization formula for compact Lie groups $G$ is very well known,
but the problem with previously existed proofs is that they rely
heavily on compactness of the acting group $G$ and none of them
generalizes to a noncompact group setting.
In the future, the author intends to adapt this proof and obtain
a localization formula which can be applied to some nontrivial
integrals that cannot be reduced to the compact group case.

\separate

The proof which appears here is a significant modification of the
localization argument which appeared in my Ph.D. thesis \cite{L1}.
This thesis contains a result on localizations
of actions of noncompact groups, which is probably the only existing one
on this subject.
Article \cite{L2} gives a very accessible introduction to \cite{L1}
and explains key ideas used there by way of examples.

\separate

I would like to thank Sergey Vasilyev and Michele Vergne for pointing
out some very serious errors in an earlier version of this article.

\separate

\begin{section}
{Equivariant Forms}
\end{section}

In this article we use the same notations as in \cite{BGV}.

\separate

Let $M$ be a ${\cal C}^{\infty}$-manifold of dimension $n$ with
an action of a (possibly noncompact) Lie group $G$,
and let $\g g$ be the Lie algebra of $G$.
The group $G$ acts on ${\cal C}^{\infty}(M)$ by the formula
$(g \cdot \phi) (x) = \phi(g^{-1} x)$.
For $X \in \g g$, we denote by $X_M$ the vector field on $M$ given by
(notice the minus sign)
$$
(X_M \cdot \phi) (x) =
\frac d{d\epsilon} \phi \bigl( \exp (-\epsilon X)x \bigr) \Bigr|_{\epsilon=0}.
$$

Let ${\cal A}(M)$ denote the (graded) algebra of smooth differential
forms on $M$, and let ${\cal C}^{\infty}(\g g) \hat \otimes {\cal A}(M)$
denote the algebra of all smooth ${\cal A}(M)$-valued functions on $\g g$.
The group $G$ acts on an element
$\alpha \in {\cal C}^{\infty}(\g g) \hat \otimes {\cal A}(M)$
by the formula
$$
(g \cdot \alpha)(X) = g \cdot (\alpha ( g^{-1} \cdot X))
\qquad \text{for all $g \in G$ and $X \in \g g$.}
$$
Let ${\cal A}^{\infty}_G(M) =
({\cal C}^{\infty}(\g g) \hat \otimes {\cal A}(M))^G$
be the subalgebra of $G$-invariant elements.
An element $\alpha$ of ${\cal A}^{\infty}_G(M)$ satisfies the relation
$ \alpha(g \cdot X) = g \cdot \alpha(X)$ and is called and
{\em equivariant differential form}.

We define the equivariant exterior differential $d_{\g g}$ on
${\cal C}^{\infty}(\g g) \hat \otimes {\cal A}(M)$ by the formula
$$
(d_{\g g} \alpha) (X) = d(\alpha(X)) - \iota(X_M) (\alpha(X)),
$$
where $\iota(X_M)$ denotes contraction by the vector field $X_M$.
This differential $d_{\g g}$ preserves ${\cal A}^{\infty}_G(M)$,
and $(d_{\g g})^2 \alpha =0$ for all $\alpha \in {\cal A}^{\infty}_G(M)$.
The elements of ${\cal A}^{\infty}_G(M)$ such that
$d_{\g g} \alpha =0$ are called {\em equivariantly closed forms}.

\separate

\begin{ex} \label{sigma}
{\em Let $T^*M$ be the cotangent bundle of $M$, and
let $\sigma$ denote the canonical symplectic form on $T^*M$.
It is defined, for example, in \cite{KaScha}, Appendix A2.
The action of the Lie group $G$ on $M$ naturally extends
to $T^*M$. Then we always have a canonical equivariantly closed
form on $T^*M$, namely, $\mu + \sigma$.
Here $\mu: \g g \to {\cal C}^{\infty}(T^*M)$ is the moment map
defined by:
$$
\mu(X): \xi \mapsto \langle \xi, X_M \rangle,
\qquad X \in \g g,\: \xi \in T^*M. \qquad \qquad \square
$$}
\end{ex}

If $\alpha$ is a non-homogeneous equivariant differential form,
$\alpha_{[k]}$ denotes the homogeneous component of degree $k$.
If $M$ is a compact oriented manifold, we can integrate equivariant
differential forms over $M$, obtaining a map
$$
\int_M : {\cal A}^{\infty}_G(M) \to {\cal C}^{\infty}(\g g)^G,
$$
defined by the formula $(\int_M \alpha) (X) = \int_M \alpha(X)_{[n]}$,
where $n=\dim M$.

\separate

Finally, notice that if $\alpha$ is an equivariantly closed form
whose top homogeneous component has degree $k$, then
$\alpha(X)_{[k]}$ is closed with respect to the ordinary
exterior differential.

\separate

\begin{section}
{Localization Formula}
\end{section}

From now on we assume that the Lie group $G$ is {\em compact}.

We recall some more notations from \cite{BGV}.
Let $M_0(X)$ be the set of zeroes of the vector field $X_M$.
We will state and prove the localization formula in the important
special case where $X_M$ has isolated zeroes.
Here, at each point $p \in M_0(X)$, the infinitesimal action of $X$ on
${\cal C}^{\infty} (M)$ gives rise to a linear transformation
$L_p$ on $T_pM$.

Because the Lie group $G$ is compact, the transformation $L_p$
is invertible and has only imaginary eigenvalues.
Thus the dimension of $M$ is even and there exists an oriented basis
$\{e_1,\dots, e_n\}$ of $T_pM$ such that for $1\le i \le l = n/2$,
$$
L_p e_{2i-1} = \lambda_{p,i} e_{2i}, \quad
L_p e_{2i} = -\lambda_{p,i} e_{2i-1}.
$$
We have $\det(L_p)=\lambda_{p,1}^2 \lambda_{p,2}^2 \dots \lambda_{p,l}^2$,
and it is natural to take the following square root (dependent only on the
orientation of the manifold):
$$
{\det}^{1/2}(L_p) = \lambda_{p,1} \dots \lambda_{p,l}.
$$

\separate

For convenience, we restate Theorem 7.11 from \cite{BGV}.

\begin{thm}  \label{BGV}
Let $G$ be a compact Lie group with Lie algebra $\g g$
acting on a compact oriented manifold $M$, and let $\alpha$
be an equivariantly closed differential form on $M$.
Let $X \in \g g$ be such that the vector field $X_M$ has
only isolated zeroes. Then
$$
\int_M \alpha(X) = (-2\pi)^l
\sum_{p \in M_0(X)} \frac {\alpha(X)(p)}{\det^{1/2} (L_p)},
$$ where $l= \dim(M)/2$, and by $\alpha(X)(p)$ we mean the value of
the function $\alpha(X)_{[0]}$ at the point $p \in M$.
\end{thm}

\separate

\begin{section}
{Proof of the Localization Formula}
\end{section}

\noindent {\it Proof of Theorem \ref{BGV}.}
Let $\pi: T^*M \twoheadrightarrow M$ denote the projection map.
We regard $M$ as a submanifold of $T^*M$ via the zero section inclusion.
Recall the canonical equivariantly closed form $\mu + \sigma$ on $T^*M$
constructed in Example \ref{sigma}. We can consider the form
\begin{equation}  \label{omega}
\omega(X) = e^{\mu(X) + \sigma} \wedge \pi^* \bigl( \alpha(X) \bigr),
\qquad X \in \g g.
\end{equation}
It is an equivariantly closed form on $T^*M$ for the reason that it is
``assembled'' from equivariantly closed forms.
Moreover, its restriction to $M$ is just $\alpha(X)$.
We will see later that, in a way, $\omega$ is the most natural
equivariant extension of $\alpha$ to $T^*M$.

\separate

Fix a $G$-invariant Riemannian metric on $T^*M$ and let $\|.\|$ be the norm.
It induces a vector bundle diffeomorphism $TM \to T^*M$ and a norm
on $TM$ which we also denote by $\|.\|$.
Let $s: M \to T^*M$ denote the image of the vector field
$-X_M$ under this diffeomorphism. Observe that
$$
\mu(X)(s(m)) = - \|X_M\|^2 \le 0, \qquad \text{for all $m \in M$},
$$
and that the set
$$
\{ m \in M;\, \mu(X)(s(m))=0 \} = \{ m \in M;\, X_M=0 \} = M_0(X)
$$
is precisely the set of zeroes of $X_M$.

\separate

Pick an $R>0$ and let $t$ increase from 0 to $R$.
Then the family of sections $\{ts\}_{t\in[0,R]}$
provides a deformation of the initial cycle $M$ into $Rs$.
In other words, there exists a chain $C_R$ in $T^*M$
of dimension $(\dim M +1)$ such that
\begin{equation}  \label{deformation}
\partial C_R = M - Rs.
\end{equation}
Note that the support of $C_R$, $|C_R|$, lies inside the set
\begin{equation}  \label{<0}
\{ \xi \in T^*M;\, \mu(X)(\xi) = \langle \xi, X_M \rangle \le 0 \}.
\end{equation}
This will ensure good behavior of the integrand when we let $R$
tend to infinity.

Observe that if we let $R \to \infty$, then we obtain a conic Borel-Moore
chain $C_{\infty}$ such that
$$
\partial C_{\infty} = M - \sum_{p \in M_0(X)} T^*_pM,
$$
where each cotangent space $T^*_pM$ is provided with
appropriate orientation.
But we will not use this observation in the proof because one cannot
interchange the order of integration and taking limit as $R \to \infty$.

\separate

\begin{lm}  \label{zero}
$\bigl( d \omega(X)_{[n]} \bigr) |_{|C_R|} =0$.
\end{lm}

\begin{rem} {\em
The form $\omega(X)_{[n]}$ itself need not be closed; this lemma
only says that $\omega(X)_{[n]}$ becomes closed when restricted to
$|C_R|$.

This lemma is the key moment in the proof and it is the only place
where we use that the form $\alpha$ is equivariantly closed and that
our Riemannian metric is $G$-invariant. Its proof is quite technical.
}\end{rem}

\pf
The form $\omega(X)_{[n]}$ can be written more explicitly as
\begin{multline*}
\omega(X)_{[n]}
= \bigl( e^{\mu(X) + \sigma} \wedge \alpha(X) \bigr)_{[n]}  \\
= e^{\mu(X)} \bigl( (e^{\sigma})_{[0]} \wedge \alpha(X)_{[n]} +
(e^{\sigma})_{[2]} \wedge \alpha(X)_{[n-2]} +
(e^{\sigma})_{[4]} \wedge \alpha(X)_{[n-4]} +\dots \bigr) \\
= e^{\mu(X)} \Bigl(  \frac 1{0!} \sigma^{0} \wedge \alpha(X)_{[n]} +
\frac 1{1!} \sigma^1 \wedge \alpha(X)_{[n-2]}+
\frac 1{2!} \sigma^2 \wedge \alpha(X)_{[n-4]}+ \dots \Bigr)  \\
= e^{\mu(X)} \sum_{k=0}^l \frac 1{k!} \sigma^k \wedge \alpha(X)_{[n-2k]}.
\end{multline*}

It is sufficient to prove that the form $\omega(X)_{[n]} |_{|C_R|}$
is closed away from zeroes of $X_M$, i.e. on the set $|C_R| \setminus M_0(X)$.
So let us pick $\xi \in |C_R| \setminus M_0(X)$,
say $\xi \in T^*_qM$, $q \notin M_0(X)$.
Because $X_M(q) \ne 0$ we can trivialize vector field $X_M$ near $q$.
That is, there is a neighborhood $U_q$ of $M$ containing
point $q$ and positively oriented local coordinates
$(x_{q,1},\dots, x_{q,n})$ defined on $U_q$ such that $X_M$ is expressed by
\begin{equation}  \label{X}
X_M = \frac {\partial}{\partial x_{q,1}}.
\end{equation}
Define a function $r: |C_R| \to \BB R$ by
$r(c) = \langle c, \mu(X) \rangle$, $c \in |C_R| \subset T^*M$.
Then $(x_{q,1},\dots, x_{q,n},r)$ defined on $\pi^{-1}(U_q) \cap |C_R|$
is a system of coordinates on $|C_R|$ containing $\xi$.

We expand the coordinate system $(x_{q,1},\dots, x_{q,n})$ of $M$
defined on $U_q$ to the standard coordinate system 
$(x_{q,1},\dots, x_{q,n},\zeta_{q,1},\dots,\zeta_{q,n})$ of $T^*M$
defined on $T^*U_q$ so that every element of $T^*U_q$ is
expressed in these coordinates as
$(x_{q,1},\dots, x_{q,n},\zeta_{q,1} dx_{q,1}+\dots+\zeta_{q,n} dx_{q,n})$.
In these coordinates,
$$
\sigma = dx_{q,1} \wedge d\zeta_{q,1} + \dots + dx_{q,n} \wedge d\zeta_{q,n}.
$$

Using this coordinate system
$(x_{q,1},\dots, x_{q,n},\zeta_{q,1},\dots,\zeta_{q,n})$
we can regard $X_M$ as a horizontal vector field on $T^*U_q$.
Because we chose our Riemannian metric to be $G$-invariant,
the vector field $X_M$ is tangent to $|C_R|$.

Using our assumption that the form $\alpha$ is equivariantly closed
we obtain:
\begin{multline*}
d \bigl( \omega(X)_{[n]} |_{|C_R|} \bigr)
= d \Bigl(
e^r \sum_{k=0}^l \frac 1{k!} \sigma^k \wedge
\alpha(X)_{[n-2k]} \Bigr) \Bigr|_{|C_R|}  \\
= dr \wedge \Bigl( e^r
\sum_{k=0}^l \frac 1{k!} \sigma^k \wedge \alpha(X)_{[n-2k]}
\Bigr) \Bigr|_{|C_R|}
+e^r d \Bigl( \sum_{k=0}^l \frac 1{k!} \sigma^k \wedge
\alpha(X)_{[n-2k]} \Bigr) \Bigr|_{|C_R|}  \\
= e^r \Bigl( \sum_{k=0}^l
\frac 1{k!} dr \wedge \sigma^k \wedge \alpha(X)_{[n-2k]}
+ \sum_{k=0}^l \frac 1{k!} \sigma^k \wedge
d\alpha(X)_{[n-2k]} \Bigr) \Bigr|_{|C_R|}  \\
= e^r \Bigl( \sum_{k=0}^l
\frac 1{k!} dr \wedge \sigma^k \wedge \alpha(X)_{[n-2k]}
+ \sum_{k=1}^l
\frac 1{k!} \sigma^k \wedge \iota(X_M) \bigl( \alpha(X)_{[n-2k+2]} \bigr)
\Bigr) \Bigr|_{|C_R|}  \\
= e^r \Bigl( \sum_{k=0}^l
\frac 1{k!} dr \wedge \sigma^k \wedge \alpha(X)_{[n-2k]}
- \sum_{k=1}^{l+1}
\frac 1{k!} \iota(X_M) \bigl(\sigma^{k} \bigr) \wedge \alpha(X)_{[n-2k+2]}
\Bigr) \Bigr|_{|C_R|}
\end{multline*}
because $\sigma^{l+1}=0$ and
$\bigl( \sigma^{k} \wedge \alpha(X)_{[n-2k+2]} \bigr) \bigr|_{|C_R|} =0 $
for dimension reasons, and so
\begin{multline*}
0 = \iota(X_M)
\bigl( \sigma^{k} \wedge \alpha(X)_{[n-2k+2]} \bigr) \bigr|_{|C_R|}  \\
= \iota(X_M)
\bigl( \sigma^{k} \bigr) \wedge \alpha(X)_{[n-2k+2]} \bigr|_{|C_R|}
+ \sigma^{k} \wedge  \iota(X_M) \bigl( \alpha(X)_{[n-2k+2]} \bigr)
\bigr|_{|C_R|}.
\end{multline*}
Hence:
\begin{multline*}
d \bigl( \omega(X)_{[n]} |_{|C_R|} \bigr)
= e^r \Bigl( \sum_{k=0}^l
\frac 1{k!} dr \wedge \sigma^k \wedge \alpha(X)_{[n-2k]}  \\
- \sum_{k=1}^{l+1}
\frac 1{(k-1)!} \iota(X_M) \bigl(\sigma \bigr)
\wedge \sigma^{k-1} \wedge \alpha(X)_{[n-2k+2]} \Bigr) \Bigr|_{|C_R|}  \\
= e^r \Bigl( \sum_{k=0}^l
\frac 1{k!} dr \wedge \sigma^k \wedge \alpha(X)_{[n-2k]}
- \sum_{k=0}^l \frac 1{k!} d\zeta_{q,1}
\wedge \sigma^{k} \wedge \alpha(X)_{[n-2k]} \Bigr) \Bigr|_{|C_R|}  \\
= e^r \Bigl( \sum_{k=0}^l
\frac 1{k!} (dr - d\zeta_{q,1}) \wedge \sigma^k \wedge \alpha(X)_{[n-2k]}
\Bigr) \Bigr|_{|C_R|} =0
\end{multline*}
because equation (\ref{X}) implies that in these coordinates,
for $c \in |C_R|$,
$$
\zeta_{q,1}(c) = \langle c, X_M \rangle = \langle c, \mu(X) \rangle = r(c),
$$
hence $dr = d\zeta_{q,1}$.
This finishes the proof of Lemma \ref{zero}.
\qed

Therefore,
\begin{equation} \label{R}
\int_M \alpha(X) = \int_M \omega(X) = \int_{Rs} \omega(X).
\end{equation}

\separate

For each $p \in M_0(X)$, there exists a positively oriented orthonormal basis
$\{e_1,\dots, e_n\}$ of $T_pM$ such that for $1\le i \le l = n/2$,
$$
L_p e_{2i-1} = \lambda_{p,i} e_{2i}, \quad
L_p e_{2i} = -\lambda_{p,i} e_{2i-1}.
$$
Using the exponential map we can find a neighborhood $U_p$ of $M$ containing
point $p$ and local coordinates $(x_{p,1},\dots, x_{p,n})$
defined on $U_p$ and centered at $p$ such that
$$
\frac{\partial}{\partial x_{p,1}}(p) = e_1,\dots,
\frac{\partial}{\partial x_{p,n}}(p) = e_n
$$
(in particular, the coordinate system $(x_{p,1},\dots, x_{p,n})$ is
positively oriented) and the vector field $X_M$ is expressed by
\begin{multline*}
X_M =
\lambda_{p,1} \Bigl( x_{p,2} \frac {\partial}{\partial x_{p,1}}
- x_{p,1} \frac {\partial}{\partial x_{p,2}}
\Bigr) + \dots \\
+ \lambda_{p,l} \Bigl( x_{p,2l} \frac {\partial}{\partial x_{p,2l-1}}
- x_{p,2l-1} \frac {\partial}{\partial x_{p,2l}} \Bigr).
\end{multline*}
Since $M_0(X)$ is discrete, we can make these neighborhoods
$U_p$'s small enough so that they do not overlap.

\separate

Find a sufficiently small $\epsilon >0$ such that,
for each $p \in M_0(X)$, the cube
$$
\{ (t_1,\dots,t_n);\, |t_1| \le \epsilon, \dots, |t_n| \le \epsilon \}
$$
is contained in the image of $(x_{p,1},\dots, x_{p,n})$.
For each $p \in M_0(X)$, we define a subset of $M$ containing $p$
$$
I_p = (x_{p,1},\dots, x_{p,n})^{-1} \bigl(
\{ (t_1,\dots,t_n);\, |t_1| \le \epsilon, \dots, |t_n| \le \epsilon \} \bigr).
$$

In formula (\ref{R}) we let $R$ tend to infinity.
Recall the definition of the equivariant form $\omega$ (formula (\ref{omega})).
Because of the presence of the term $e^{\mu(X)+\sigma}$ there and
because of the exponential decay on the set (\ref{<0}) we get:
\begin{multline*}
\int_M \alpha(X) = \lim_{R \to \infty} \int_{Rs} \omega(X)
= \lim_{R \to \infty}
\int_{Rs \cap \bigcup_{p \in M_0(X)} \pi^{-1} I_p} \omega(X)  \\
= \sum_{p \in M_0(X)} \lim_{R \to \infty}
\int_{Rs \cap \pi^{-1} I_p} \omega(X).
\end{multline*}
That is, when we integrate over $Rs$, only the covectors whose basepoint
lies in some $I_p$ (i.e. ``near'' $M_0(X)$) count.
Thus we see that the integral is localized around the
set of zeroes $M_0(X)$.
It remains to calculate the contribution of each point $p$:
$$
\lim_{R \to \infty} \int_{Rs \cap \pi^{-1} I_p} \omega(X).
$$

\separate

Recall that
$$
\Bigl\{ \frac {\partial}{\partial x_{p,1}} = e_1, \dots,
\frac {\partial}{\partial x_{p,n}} =e_n \Bigr\}
$$
is an orthonormal basis of $T^*_pM$.
Hence the Riemannian metric around each $p$ can be
written in terms of the frame on $U_p$
$$
\Bigl( \frac {\partial}{\partial x_{p,1}}, \dots,
\frac {\partial}{\partial x_{p,n}} \Bigr)
$$
as $I + E_p(m)$, where $I$ is the $n \times n$ identity matrix
and $E_p(m)$ is a symmetric $n \times n$ matrix which depends on $m$
and vanishes at the point $p$.

We expand each coordinate system $(x_{p,1},\dots, x_{p,n})$ of $M$
defined on $U_p$ to the standard coordinate system 
$(x_{p,1},\dots, x_{p,n},\zeta_{p,1},\dots,\zeta_{p,n})$ of $T^*M$
defined on $T^*U_p$ so that every element of $T^*U_p$ is
expressed in these coordinates as
$(x_{p,1},\dots, x_{p,n},\zeta_{p,1} dx_{p,1}+\dots+\zeta_{p,n} dx_{p,n})$.

In these coordinates,
\begin{multline*}
s = \bigl( x_{p,1},\dots, x_{p,n},
\lambda_{p,1}(- x_{p,2} dx_{p,1} + x_{p,1} dx_{p,2}) + \dots  \\
+ \lambda_{p,l} (- x_{p,2l} dx_{p,2l-1} + x_{p,2l-1} dx_{p,2l})
+ O(x^2) \bigr),
\end{multline*}
where $O(x^2)$ denotes terms of quadratic order in $x_{p,1},\dots, x_{p,n}$.
We also have:
$$
\sigma = dx_{p,1} \wedge d\zeta_{p,1}+ \dots + dx_{p,n} \wedge d\zeta_{p,n},
$$
and
\begin{multline*}
\mu(X) : (x_{p,1},\dots, x_{p,n},\zeta_{p,1},\dots,\zeta_{p,n})  \\
\mapsto
\lambda_{p,1}(x_{p,2} \zeta_{p,1} -x_{p,1} \zeta_{p,2}) + \dots
+ \lambda_{p,l}( x_{p,2l} \zeta_{p,2l-1} - x_{p,2l-1} \zeta_{p,2l}).
\end{multline*}

\separate

Let us choose a coordinate system for $Rs \cap \pi^{-1} I_p$, namely
$$
(y_{p,1}=x_{p_1} \circ \pi, \dots, y_{p,n} = x_{p,n} \circ \pi).
$$
It is positively oriented. Then in these coordinates
$$
dx_{p,1}=dy_{p,1}, \dots, dx_{p,n}=dy_{p,n},
$$
\begin{multline*}
d\zeta_{p,1} = - R\lambda_{p,1} dy_{p,2} + R \cdot O(y),  \qquad
d\zeta_{p,2} = R\lambda_{p,1}dy_{p,1} + R \cdot O(y), \dots,  \\
d\zeta_{p,2l-1} = - R\lambda_{p,l}dy_{p,2l} + R \cdot O(y), \quad
d\zeta_{p,2l} = R\lambda_{p,l}dy_{p,2l-1} + R \cdot O(y),
\end{multline*}
$$
\mu(X) |_{Rs} = -R\lambda_{p,1}^2 (y_{p,1}^2 + y_{p,2}^2) - \dots
- R\lambda_{p,l}^2 (y_{p,2l-1}^2 + y_{p,2l}^2)
+ R \cdot O(y^3).
$$
Thus $\sigma|_{Rs}$ in these coordinates becomes
$$
\sigma|_{Rs} = -2R \lambda_{p,1} dy_{p,1} \wedge dy_{p,2} - \dots
-2R \lambda_{p,l} dy_{p,2l-1} \wedge dy_{p,2l} + R \cdot O(y).
$$
Hence
$$
(e^{\sigma})_{[n]} |_{Rs} =
(-2R)^l \lambda_{p,1} \dots \lambda_{p,l} dy_{p,1} \wedge \dots \wedge dy_{p,n}
+ R^l \cdot O(y).
$$

On the other hand, the expression of $\pi^* \alpha(X) |_{Rs}$ in these
coordinates does not depend on $R$.
We have:
\begin{multline*}
\omega(X)_{[n]}  \\
= \bigl( e^{\mu(X)+\sigma} \bigr)_{[n]} \alpha(X)_{[0]} +
\bigl( e^{\mu(X)+\sigma} \bigr)_{[n-2]} \alpha(X)_{[2]}+ \dots +
\bigl( e^{\mu(X)+\sigma} \bigr)_{[0]} \alpha(X)_{[n]}.
\end{multline*}
Hence
\begin{multline*}
\omega(X)_{[n]} |_{Rs}
= \bigl( e^{\mu(X)+\sigma} \bigr)_{[n]} \alpha(X)_{[0]} |_{Rs}  \\
+ e^{\mu(X)}
\bigl( \text{terms depending on $R$ through powers of $R$
strictly less than $\frac n2$} \bigr)  \\
= e^{\mu(X)} (e^{\sigma})_{[n]}
\alpha(X)_{[0]} |_{Rs} + e^{\mu(X)}o(R^{n/2}) = \\
e^{- R \bigl( \lambda_{p,1}^2 (y_{p,1}^2 + y_{p,2}^2) + \dots
+ \lambda_{p,l}^2 (y_{p,2l-1}^2 + y_{p,2l}^2) + O(y^3) \bigr)}
\cdot \bigl( (e^{\sigma})_{[n]} \alpha(X)_{[0]} |_{Rs} + o(R^{n/2}) \bigr).
\end{multline*}

\separate

Recall that the integration takes place over $I_p$,
hence the range of variables is
$|y_{p,1}| \le \epsilon, \dots, |y_{p,n}| \le \epsilon$.
We perform a change of coordinates
$$
z_{p,1} = \sqrt{R} y_{p,1}, \dots, z_{p,n} = \sqrt{R} y_{p,n}, \qquad
|z_{p,1}| \le \epsilon \sqrt{R}, \dots, |z_{p,n}| \le \epsilon\sqrt{R}.
$$
In these new coordinates,
\begin{multline*}
\omega(X)_{[n]} |_{Rs} =  
e^{- \bigl( \lambda_{p,1}^2 (z_{p,1}^2 + z_{p,2}^2) + \dots
+ \lambda_{p,l}^2 (z_{p,2l-1}^2 + z_{p,2l}^2) + O(z^3)/O(\sqrt{R}) \bigr)}  \\
\cdot \bigl( (e^{\sigma})_{[n]} \alpha(X)_{[0]} |_{Rs} + 1/O(\sqrt{R}) \bigr),
\end{multline*}
and
$$
(e^{\sigma})_{[n]} |_{Rs} =
(-2)^l \lambda_{p,1} \dots \lambda_{p,l} dz_{p,1} \wedge \dots \wedge dz_{p,n}
+ O(z)/O(\sqrt{R}).
$$

As $R$ tends to infinity, the integrand tends pointwise to
\begin{multline*}
(-2)^l \lambda_{p,1} \dots \lambda_{p,l}
e^{- \lambda_{p,1}^2 (z_{p,1}^2 + z_{p,2}^2) - \dots
- \lambda_{p,l}^2 (z_{p,n-1}^2 + z_{p,n}^2)}  \\
\cdot \alpha(X)_{[0]}(p) dz_{p,1} \wedge \dots \wedge dz_{p,n}.
\end{multline*}
Hence by the Lebesgue dominated convergence theorem,
\begin{multline*}
\int_M \alpha(X) =
\sum_{p \in M_0(X)} \lim_{R \to \infty}
\int_{Rs \cap \pi^{-1} I_p} \omega(X)  \\
= \sum_{p \in M_0(X)}
\int_{\BB R^n} (-2)^l \lambda_{p,1} \dots \lambda_{p,l} \alpha(X)_{[0]}(p)  \\
\cdot e^{- \lambda_{p,1}^2 (z_{p,1}^2 + z_{p,2}^2) - \dots
- \lambda_{p,l}^2 (z_{p,n-1}^2 + z_{p,n}^2)}
dz_{p,1} \wedge \dots \wedge dz_{p,n}  \\
= \sum_{p \in M_0(X)}
(-2\pi)^l \frac 1{\lambda_{p,1} \dots \lambda_{p,l}} \alpha(X)_{[0]}(p)
= (-2\pi)^l \sum_{p \in M_0(X)}
\frac {\alpha(X)(p)}{{\det}^{1/2}(L_p)}.
\end{multline*}
This finishes our proof of Theorem \ref{BGV}.
\qed

\separate

\noindent
{matvei@math.umass.edu}

\noindent
{DEPARTMENT OF MATHEMATICS AND STATISTICS, UNIVERSITY OF MASSACHUSETTS,
LEDERLE GRADUATE RESEARCH TOWER, 710 NORTH PLEASANT STREET, AMHERST,
MA 01003}

\end{document}